\documentstyle[twoside,leqno,11pt]{article}

\newfont{\Bbb}{msbm10 at 11pt}
\newfont{\Bbm}{msbm10 at 8pt}
\newfont{\Bbs}{msbm10 at 7pt}

\newtheorem{thm}{Theorem}[section]
\newtheorem{lem}[thm]{Lemma}

\begin{document}

\title{\bf A Floquet-Liapunov theorem in Fr\'echet spaces
         }

\author{George N. Galanis - Efstathios E. Vassiliou}
\date{}

\maketitle

\begin{abstract}
\noindent
Based on \cite{4}, we prove a variation of the theorem in title, 
for equations with periodic coefficients, in Fr\'echet spaces. The 
main result gives equivalent conditions ensuring the reduction of 
such an equation to one with constant coefficient. In the particular
case of $\mbox{\Bbb C}^{\infty}$, we obtain the exact analogue of the
classical theorem. Our approach essentially uses the fact that
a Fr\'echet space is the limit of a projective sequence of Banach 
spaces. This method can also be applied for a geometric 
interpretation of the same theorem within the context of total 
differential equations in Fr\'echet fiber bundles. 
\end{abstract}

\section*{Introduction}

In spite of the great progress in the study of topological 
vector spaces and their rich structure, a general solvability
theory of differential equations in the non Banach case is still 
missing. This led many authors (see e.g. R.S.~Hamilton \cite{5}, R.~Lemmert 
\cite{7} and the survey article \cite{2}) to study special classes of 
differential equations where some types of solutions can be found. 
Also in this context, G. Galanis in \cite{4}, published 
in this journal, has studied a class of linear differential 
equations in a Fr\'echet space $\mbox{{\Bbb F}}$, which can be always
uniquely solved under given initial conditions. The key to this approach
is the fact that $\mbox{{\Bbb F}}$ can be though of as the limit of a
projective system of Banach spaces: $\mbox{{\Bbb F}} \equiv
\displaystyle{\lim_{\longleftarrow}}\mbox{{\Bbb E}}_i$.
As a matter of fact, this class consists of all equations
$\dot{x}=A(t)\cdot x$ whose coefficient $A$ can be
factored in the form 
\[
A = \varepsilon \circ A^*,
\]
where $A^* : [0,1] \longrightarrow H(\mbox{{\Bbb F}})$ is a continuous map,
$H(\mbox{{\Bbb F}})$ is the Fr\'echet subspace of
$\prod_{i \in \mbox{{\Bbs N}}}({\cal L}(\mbox{{\Bbb E}}_i))$ containing
all the sequences $(f_i)_{i \in {\Bbb N}}$ which form projective systems and
$\varepsilon : H(\mbox{{\Bbb F}}) \longrightarrow {\cal L}(\mbox{{\Bbb F}})
: (f_i)_{i \in \mbox{{\Bbs N}}} \mapsto
\displaystyle{\lim_{\longleftarrow}}f_i$.

Though the Cauchy problem for linear equations, as above, 
could be studied yet by more advanced methods (referring to 
equations in locally convex spaces, cf. e.g. \cite{2,7}), \cite{4} 
expounds an elementary approach leading to an explicit 
description of the solutions (apart from their theoretical 
existence) as projective limits of solutions in Banach spaces, 
where the corresponding theory is fairly complete. 

A natural question that arises now is whether we can exploit 
the same approach to obtain a kind of a Floquet-Liapunov 
theorem in the context of Fr\'echet spaces. 

As it is well known, the classical theorem asserts that a 
differential equation $\dot{x} = A(t) \cdot x$, with $A(t)$ a periodic
complex matrix, can always be reduced to an equation with constant
coefficient (cf. e.g. \cite{1,10,16}). The key to this reduction is the 
existence of the logarithm of $\Phi(\tau)$, if $\Phi$  is the fundamental 
solution of the equation and $\tau$ the period. 

However, this result is not true, in general, even in the 
case of Banach spaces, since the existence of such a logarithm is 
not always ensured. In this respect we refer to \cite{9} and \cite{12}, 
where the reader can also find some sufficient conditions for the 
validity of the result in this context. 

In the present paper, which is a natural outgrowth of \cite{4}, we 
present a variation of the theorem in title. More precisely, in 
the first main result (Theorem 2.3) we obtain equivalent 
conditions implying the reduction of a periodic equation to one 
with constant coefficient. The basic idea, in this direction, is to 
consider a generalized monodromy homomorphism of the equation at 
hand, in which the pathological group $GL(\mbox{{\Bbb F}})$ (fundamentally 
involved in the classical or Banach case) is replaced by an 
appropriate subgroup $H^o(\mbox{{\Bbb F}})$ of
$H(\mbox{{\Bbb F}})$ fully described in Section 1,
since the classical homomorphism is useless in our context. 

For the sake of completeness, we apply the main Theorem 2.3 
also to the case of $\mbox{\Bbb C}^{\infty}$. It is worthy to note that,
although we still work in a Fr\'echet space, now the reduction under
discussion always holds true. Therefore (see Theorem 3.2), we obtain the 
exact analogue of the classical theory. 

The results of this note can be also profitably used for a 
generalization of the Floquet-Liapunov theory within the geometric 
context of Fr\'echet fiber bundles. In the latter, ordinary linear 
equations correspond to equations with total differentials which, 
in turn, are interpreted as flat connections. Hence, the 
generalized analogue of the main result leads to the reduction of 
a flat connection to one with ``constant" coefficients on a trivial 
bundle, as we briefly discuss in Section 4.

\section{Preliminaries}   

Let $\mbox{{\Bbb F}}$ be a Fr\'echet space. It is well known (see for
instance \cite{11}) that $\mbox{{\Bbb F}}$ can be identified with the limit
of a projective system
$\{\mbox{{\Bbb E}}_i;\rho_{ji}\}_{i,j \in \mbox{{\Bbs N}}}$ of
Banach spaces, i.e. 
$\mbox{{\Bbb F}} \equiv
\displaystyle{\lim_{\longleftarrow}}\mbox{{\Bbb E}}_i$. More precisely,
each $\mbox{{\Bbb E}}_i$ is the completion of the normed space
$\mbox{{\Bbb F}}/Ker(p_i)$, where $\{p_i\}_{i \in \mbox{{\Bbs N}}}$
are the seminorms of $\mbox{{\Bbb F}}$ and $\rho_{ji}$ are the
extensions of the natural projections
\[
\mbox{{\Bbb F}}/Ker(p_j) \longrightarrow \mbox{{\Bbb F}}/Ker(p_i) :
x+Ker(p_j) \mapsto  x+Ker(p_i), \quad j \geq i.
\]

As explained in the Introduction, in order to study linear differential
equations in a Fr\'echet space $\mbox{{\Bbb F}}$ as above, we consider
the following space:
\[
H(\mbox{{\Bbb F}}) := \{ (f_i)_{i \in \mbox{{\Bbs N}}} :
f_i \in {\cal L}(\mbox{{\Bbb E}}_i)\; \mbox{ and }\;
\displaystyle{\lim_{\longleftarrow}}f_i\;
\mbox{ exists}\}.
\]
It can be proved that $H(\mbox{{\Bbb F}}) \equiv
\displaystyle{\lim_{\longleftarrow}}(H_k(\mbox{{\Bbb F}}))$,
where each
\[
H_k(\mbox{{\Bbb F}}) := \{ (f_1,...,f_k) : f_i \in
{\cal L}(\mbox{{\Bbb E}}_i) \,\mbox{ and }\, \rho_{ji} \circ f_j =
f_i \circ \rho_{ji} \;  (j \geq i) \} ;
\; k \in \mbox{{\Bbb N}},
\]
is a Banach space. Thus $H(\mbox{{\Bbb F}})$ is a Fr\'echet space. 

Similarly, we set 
\[
H^o(\mbox{{\Bbb F}}) := H(\mbox{{\Bbb F}}) \cap
\prod_{i \in \mbox{{\Bbs N}}} GL(\mbox{{\Bbb E}}_i).
\]
This is is a topological group, as the limit (within an isomorphism) 
of the projective system of the Banach-Lie groups 
\[
H^o_k(\mbox{{\Bbb F}}) := H_k(\mbox{{\Bbb F}}) \cap
\prod_{i=1}^k GL(\mbox{{\Bbb E}}_i), \quad  k \in \mbox{{\Bbb N}}.
\]

In the study of differential equations the exponential 
mapping plays a fundamental role (cf. e.g. \cite{3}). However, in the 
case of a Fr\'echet space $\mbox{{\Bbb F}}$, an exponential map of
the form $\exp : {\cal L}(\mbox{{\Bbb F}}) \longrightarrow
GL(\mbox{{\Bbb F}})$ (in analogy to the Banach case) cannot be
defined, as a consequence of the incompleteness of
${\cal L}(\mbox{{\Bbb F}})$ and the pathological structure of
$GL(\mbox{{\Bbb F}})$. This difficulty is overcome by inducing a
{\em generalized exponential\/} $\mbox{Exp} : H(\mbox{{\Bbb F}})
\longrightarrow H^o(\mbox{{\Bbb F}})$ defined in the
following way: If $\mbox{{\Bbb F}} \equiv
\displaystyle{\lim_{\longleftarrow}} \mbox{{\Bbb E}}_i$, we denote by
$\exp_i : {\cal L}(\mbox{{\Bbb E}}_i) \longrightarrow
GL(\mbox{{\Bbb E}}_i)$ the ordinary exponential mapping of the
Banach-Lie group $GL(\mbox{{\Bbb E}}_i)$. We check that the finite products 
\[
\exp_1  \times  \dots  \times \exp_k : H_k(\mbox{{\Bbb F}})
\longrightarrow H_k^o(\mbox{{\Bbb F}}); \quad k \in \mbox{{\Bbb N}},
\]
form a projective system, thus we may define
\[
\mbox{Exp} := \displaystyle{\lim_{\longleftarrow}}(\exp_1 \times
\dots  \times \exp_k).
\]

Using now the differentiability of mappings between Fr\'echet 
spaces in the sense of J. Leslie (\cite{8}), the previous considerations
lead to the following main result.

\begin{thm}       
Let $\mbox{{\Bbb F}}$ be a Fr\'echet space and the (homogeneous) 
linear equation
\begin{equation}       
\dot{x} = A(t) \cdot x,
\end{equation}
where the coefficient $A : I=[0,1] \longrightarrow
{\cal L}(\mbox{{\Bbb F}})$ is continuous and factors in the form
$A=\varepsilon \circ A^*$, with $A^*:I \longrightarrow
H(\mbox{{\Bbb F}})$ being a continuous mapping and $\varepsilon :
H(\mbox{{\Bbb F}}) \longrightarrow {\cal L}(\mbox{{\Bbb F}})$ given by
$\varepsilon ((f_i))=\displaystyle{\lim_{\longleftarrow}}f_i$. Then (1)
admits a unique solution for any given initial condition. 
\end{thm}

The complete study of equations of type (1) has been already 
given in \cite{4}. However, since we need the exact form of the 
solutions, we outline the basic steps of the proof. 

\medskip
\noindent
{\bf Proof.} Let $(t_0,x_0) \in \mbox{{\Bbb R}} \times \mbox{{\Bbb F}}$
be an initial condition. Considering a projective system
$\{\mbox{{\Bbb E}}_i;\rho_{ji}\}_{i,j \in \mbox{{\Bbs N}}}$ of Banach
spaces with $\mbox{{\Bbb F}}\equiv
\displaystyle{\lim_{\longleftarrow}}\mbox{{\Bbb E}}_i$, we
obtain the equation (on $\mbox{{\Bbb E}}_i$):

\vspace{.3cm}
\noindent        
$(1.i)$ \hfill \hspace{-26pt} \centerline{$\dot{x}_i = A_i(t)\cdot x_i;
                    \quad  i \in \mbox{{\Bbb N}}$,}

\medskip
\noindent
where $A_i$ is the projection of $A$ to the $i$-factor. Since
$\mbox{{\Bbb E}}_i$ is Banach, there exists a unique solution $f_i$ of
$(1.i)$ through $(t_0,\rho_i(x_0))$, if $\rho_i : \mbox{{\Bbb F}}
\longrightarrow \mbox{{\Bbb E}}_i$ denotes the canonical projection.
For any $j \geq i$, we check that $\rho_{ji} \circ f_j$ is also a
solution of $(1.i)$ with the same initial data, hence
$\rho_{ji} \circ f_j =f_i$. As a result, the $C^1$-mapping
$f:=\displaystyle{\lim_{\longleftarrow}}f_i$ can be defined and yields
\[
\dot{f}(t) = (\dot{f}_i(t))_{i \in \mbox{{\Bbs N}}} =
(A_i(t)(f_i(t)))_{i \in \mbox{{\Bbs N}}} = A(t)(f(t)),
\]
\[
f(t_0) = (\rho_i(x_0))_{i \in \mbox{{\Bbs N}}} = x_0.
\]
Therefore, $f$ is the desired solution of (1). Moreover, $f$ is unique 
with respect to $(t_0,x_0)$, since for any other solution $g$ we may 
check, using analogous techniques as before, that $\rho_i\circ g =
\rho_i \circ f$, $i \in \mbox{{\Bbb N}}$. \hfill $\Box$
                                                  
\medskip
\noindent
{\bf Remarks.} From the preceding proof it is clear that the method 
applied therein provides explicit solutions of equation (1). 
Moreover, the present approach, apart from being elementary (in 
contrast to more general methods known for equations in arbitrary 
locally convex spaces; cf. \cite{2,7}), is very convenient for the 
geometric generalization discussed in Section 4.

\section{Floquet-Liapunov theory} 

In this section we consider again a differential equation of
the form (1) with a {\em periodic coefficient} $A$, satisfying
$A = \varepsilon \circ A^*$ as in Theorem 1.1. Without loss of
generality, we may assume that {\em the period is 1}.

\begin{lem}   
The coefficient $A$ is periodic if and only if $A^*$ is.
\end{lem}

\medskip
\noindent
{\bf Proof.} The periodicity of $A^*$ clearly implies that of $A$. 
Conversely, equality $A(t+1)=A(t)$ implies that
$\displaystyle{\lim_{\longleftarrow}}(A_i(t+1)) =
\displaystyle{\lim_{\longleftarrow}}(A_i(t))$, thus
\begin{equation}  
A_i(t+1) \circ \rho_i = A_i(t) \circ \rho_i; \quad
         i \in \mbox{{\Bbb N}}, \; t \in I,
\end{equation}
where $\rho_i : \mbox{{\Bbb F}} \longrightarrow \mbox{{\Bbb E}}_i$ are
the canonical projections. Taking into account that any element of
$\mbox{{\Bbb E}}_i$ is the limit of a sequence of elements of
$\mbox{{\Bbb F}}/Ker(p_i)$ and $\mbox{{\Bbb F}}$ is projected onto
the previous quotient (for each $i \in \mbox{{\Bbb N}}$), (2) results
in $A_i(t+1)=A_i(t)$, thus ensuring the periodicity of $A^*$.
\hfill $\Box$

\medskip
In the proof of the main Theorem we shall also need the 
following auxiliary result 

\begin{lem}   
{\em (\cite[Proposition 3.1]{4})} Let $f \in {\cal L}(\mbox{{\Bbb F}})$.
Then $f = \displaystyle{\lim_{\longleftarrow}}f_i$,
$f_i \in {\cal L}(\mbox{{\Bbb E}}_i)$, if and only if, for each
$i \in \mbox{{\Bbb N}}$, there exists $M_i > 0$ such that
\[
p_i(f(u)) \leq M_i \cdot p_i(u), \quad u \in \mbox{{\Bbb F}}.
\]
\end{lem}

From Lemma 2.1, it follows that each differential equation of
type $(1.i)$ has periodic coefficient. Therefore, the corresponding 
monodromy homomorphisms (\cite{1,6}) have the form 
\[
\alpha_i^{\#} : \mbox{{\Bbb Z}} \longrightarrow GL(\mbox{{\Bbb E}}_i)
: n \mapsto  \Phi_i(n),
\]
where $\Phi_i$  is the fundamental solution (resolvant) of $(1.i)$. 

Since the solutions of (1) are limits of solutions of the 
projective system of equations $(1.i)$, the family
$(\alpha_i^{\#}(n))_{i \in \mbox{{\Bbs N}}}$, for
each $n \in \mbox{{\Bbb Z}}$, forms a projective system. Hence
we can define the homomorphism 
\begin{equation}  
\alpha^{\#}  : \mbox{{\Bbb Z}} \longrightarrow
H^o(\mbox{{\Bbb F}}) : n \mapsto
(\alpha_i^{\#}(n))_{i \in \mbox{{\Bbs N}}}.
\end{equation}
The previous homomorphism is the key to our approach as replacing the 
classical monodromy homomorphism with values in
$GL(\mbox{{\Bbb F}})$. The latter, as is well known, is useless in
the Fr\'echet framework since it does not admit even a
topological group structure.

We are now in a position to prove the first main result of 
this note. 

\begin{thm}   
{\bf (Floquet-Liapunov)} Assume that the coefficient $A$ of 
(1) is periodic. Then the following conditions are equivalent: 

\noindent
(i) Equation (1) reduces to an (equivalent) equation
\[
\dot{y} = B \cdot y, \quad B \mbox{ constant }
\]
by means of a periodic transformation
$y = (\varepsilon \circ Q^*)(t)\cdot x$, for some continuous
$Q^* : \mbox{{\Bbb R}} \longrightarrow H^o(\mbox{{\Bbb F}})$.

\noindent
(ii) There exists $\overline{B} \in H(\mbox{{\Bbb F}})$ such that
$\mbox{Exp}(\overline{B}) = \alpha^{\#}(1)$.

\noindent
(iii) The monodromy homomorphism $\alpha^{\#}$  can be extended to a
homomorphism $F : \mbox{{\Bbb R}} \longrightarrow
H^o(\mbox{{\Bbb F}})$.
\end{thm}

\medskip
\noindent
{\bf Proof.} Assume that (i) holds true. Then, setting
$Q := \varepsilon \circ Q^*$, we see that $Q(t) =
\displaystyle{\lim_{\longleftarrow}}(Q_i(t))$, for every
$t \in \mbox{{\Bbb R}}$, where $Q_i : \mbox{{\Bbb R}}
\longrightarrow GL(\mbox{{\Bbb E}}_i)$. Since
$B \cdot Q = Q \cdot \dot{x} + \dot{Q} \cdot x$ (as a result of
the transformation $Q$) we check that
\begin{equation}     
B(u) = Q(0) (A(0) \cdot u) + \dot{Q}(0) \cdot u, \quad
u \in \mbox{{\Bbb F}}.
\end{equation}
By Lemma 2.2, along with
$A(0) = \displaystyle{\lim_{\longleftarrow}}(A_i(0))$ and
$Q(0) = \displaystyle{\lim_{\longleftarrow}}(Q_i(0))$, (4)
implies that $B = \displaystyle{\lim_{\longleftarrow}}B_i$, where
$B_i \in {\cal L}(\mbox{{\Bbb E}}_i)$.
Transforming now each $(1.i)$ via $y_i = Q_i(t) \cdot x_i$,
we observe that
\[
\dot{y}_i = \rho_i(\dot{y}) = \rho_i(B \cdot y) =
B_i \cdot \rho_i(y) = B_i \cdot y_i
\]
with $B$ constant. Therefore, applying the Theorem of Floquet in 
Banach spaces, $\exp_i(B_i) = \alpha^{\#}_i(1)$, $i \in \mbox{{\Bbb N}}$.
Hence, $\mbox{Exp}((B_i)_{i \in \mbox{{\Bbs N}}}) = \alpha^{\#}(1)$,
which proves (ii) for $\overline{B} := (B_i)_{i \in \mbox{{\Bbs N}}}$.

If (ii) holds, $\overline{B}$ has necessarily the form
$\overline{B}:=(B_i)_{i \in \mbox{{\Bbs N}}}$, with
$B_i \in {\cal L}(\mbox{{\Bbb E}}_i)$ and $\exp_i(B_i) =
\alpha_i^{\#}(1)$. Therefore, once again by Floquet's
Theorem, there exists a homomorphism $F_i : \mbox{{\Bbb R}}
\longrightarrow GL(\mbox{{\Bbb E}}_i)$ extending $\alpha_i^{\#}$.
Since $F_i(t) = \exp_i(t \cdot B_i)$, we check that
$(F_i(t))_{i \in \mbox{{\Bbs N}}}$ is a projective
system, for any $t \in \mbox{{\Bbb R}}$, thus we obtain condition
(iii) by defining
\[
F : \mbox{{\Bbb R}} \longrightarrow H^o(\mbox{{\Bbb F}}) :
t \mapsto  (F_i(t))_{i \in \mbox{{\Bbs N}}}.
\]

Finally, (iii) implies (i) as follows: For an $F$ as in the
statement, we have that $F = (F_i)_{i \in \mbox{{\Bbs N}}}$
with $F_i|{\mbox{{\Bbb Z}}} = \alpha_i^{\#}$. Thus, each equation
$(1.i)$ (with periodic coefficient) reduces to the equation with
constant coefficient
\begin{equation} 
\dot{y}_i = B_i \cdot y_i; \quad i \in \mbox{{\Bbb N}},
\end{equation}
by means of the transformation $y_i = Q_i(t) \cdot x_i$, where
\[
Q_i(t) = \exp_i(tB_i) \cdot \Phi_i^{-1}(t),
\]
and $B_i = \log_i(\alpha_i^{\#}(1))$. We check that
$(Q_i(t))_{i \in \mbox{{\Bbs N}}}$ is a projective
system, for each $t \in \mbox{{\Bbb R}}$, and we set
$Q^*(t) := (Q_i(t))_{i \in \mbox{{\Bbs N}}}$. On the other
hand, for any $j \geq i$, $B_j(u) = \dot{y}_j(0)$, where $y_j$
is the solution of (5) with initial condition $(0,u)$. Since
$x_j = Q_j(t)^{-1} \cdot y_j$ is a solution of
$(1.i)$ with coefficient $A_j(t)$, clearly $\rho_{ji} \circ x_j$
is a solution of $(1.i)$ with coefficient $A_i(t)$. Therefore,
\[
y_i = Q_i(t)(\rho_{ji} \circ x_j)
\]
is the solution of (5) with coefficient $B_i$ and initial condition 
$(0,\rho_{ji}(u))$. Moreover, $\rho_{ji} \circ y_j = y_i$. Hence,
\[
(B_i \circ \rho_{ji})(u) = \dot{y}_i(0) = \rho_{ji}(\dot{y}_j(0)) =
(\rho_{ji} \circ B_j)(u), \quad u \in \mbox{{\Bbb E}}_j.
\]
As a result, $B := \displaystyle{\lim_{\longleftarrow}}B_i  \in
{\cal L}(\mbox{{\Bbb F}})$ exists and is constant. Applying
now the transformation $y = (\varepsilon \circ Q^*)(t) \cdot x$,
we see that
\[
\dot{y}(t) = ((\displaystyle{\lim_{\longleftarrow}}Q_i(t)) \cdot
(x_i(t))_{i \in \mbox{{\Bbs N}}})^{\textstyle{.}} =
((Q_i(t)(x_i(t)))^{\textstyle{.}})_{i \in \mbox{{\Bbs N}}} =
\]
\[
= (\dot{y}_i(t))_{i \in \mbox{{\Bbs N}}} =
(B_i \cdot y_i(t))_{i \in \mbox{{\Bbs N}}} = B \cdot y(t),
\]
for every $t \in \mbox{{\Bbb R}}$, thus completing the proof.
\hfill $\Box$

\medskip
\noindent
{\bf Remark.} We note that in the previous theorem the quantities
$B$, $\overline{B}$ do not coincide, as in the Banach case, but they are
related via $\varepsilon(\overline{B}) = B$.

\section{An application: equations in $\mbox{{\Bbb C}}^{\infty}$}

In this section we shall prove that the conditions of Theorem
2.3 are always satisfied in the case of a periodic equation of type 
(1) with $\mbox{{\Bbb F}}=\mbox{{\Bbb C}}^{\infty}$. This illustrates
the elementary methods of Section 2 in the concrete case of
$\mbox{{\Bbb C}}^{\infty}$. Here $\mbox{{\Bbb C}}^{\infty}$ is viewed
as the limit $\displaystyle{\lim_{\longleftarrow}}\mbox{{\Bbb C}}^n$,
$n \in \mbox{{\Bbb N}}$, with corresponding connecting morphisms
$\rho_{ji} : \mbox{{\Bbb C}}^j \longrightarrow \mbox{{\Bbb C}}^i$
($i,j \in \mbox{{\Bbb N}}; \; j \geq i$) the natural projections. 

\smallskip
To this end we first need the following

\begin{lem}  
Let $\{ f_n : \mbox{{\Bbb C}}^n \longrightarrow
\mbox{{\Bbb C}}^n \}_{n \in \mbox{{\Bbs N}}}$ be a family of
continuous linear mappings and let $\{ M_n \in
{\cal M}_n(\mbox{{\Bbb C}}^n) \}_{n \in \mbox{{\Bbs N}}}$ be the
corresponding matrices with respect to the natural basis of
$\mbox{{\Bbb C}}^n$. Then the following conditions
are equivalent:

\noindent
(i) $\{ f_n \}_{n \in \mbox{{\Bbs N}}}$ is a projective system.

\smallskip
\noindent
(ii) \[ M_{n+1} = \left( \begin{array}{cc}
                            M_n & 0 \\
                            \mu_n & \lambda_n
                           \end{array} \right);
                 \quad  n \in \mbox{{\Bbb N}},
      \]
where $\mu_n \in \mbox{{\Bbb C}}^n$, and $\lambda_n \in \mbox{{\Bbb C}}$.
\end{lem}

\medskip
\noindent
{\bf Proof.} If (i) holds, then $\rho_{n+1,n} \circ f_{n+1} = f_n \circ
\rho_{n+1,n}$ ($n \in \mbox{{\Bbb N}}$) implies that
\[
f_{n+1}(e_k^{n+1}) = (f_n(e^n_k),\mu^n_k), \quad 1 \leq k \leq n,
\]
\[
f_{n+1}(e^{n+1}_{n+1}) = (0, \lambda_n),
\]
where $\{ e_k^n \}_{1 \leq k \leq n}$  is the natural basis of
$\mbox{{\Bbb C}}^n$ and $\mu^n_k, \, \lambda_n \in
\mbox{{\Bbb C}}$. Therefore,
\[
M_{n+1} = \left( \begin{array}{cc}
                \left( \begin{array}{c}
                f_n(e^n_k) \\
                \mu^n_k \end{array} \right)_{1 \leq k \leq n}
                & \begin{array}{c}
                0 \\
                \smallskip
                \lambda_n \end{array}
                \end{array} \right) =
\left( \begin{array}{cc}
               M_n & 0 \\ 
              (\mu^n_k)_{1 \leq k \leq n} & \lambda_n
              \end{array} \right)
\]
which gives (ii) for $\mu_n := (\mu^n_k)_{1\leq k\leq n}$.

\medskip
Conversely, if we assume (ii), then for every
$(x_1, \dots ,x_{n+1}) \in \mbox{{\Bbb C}}^{n+1}$,
\[
(\rho_{n+1,n} \circ f_{n+1})(x_1, \dots ,x_{n+1}) =
\rho_{n+1,n}(M_{n+1} \cdot \left( \begin{array}{c}
                                   x_1 \\
                                   \vdots \\
                                   x_{n+1} \end{array} \right) ) =
\]
\[
= M_n \cdot \left( \begin{array}{c}
                      x_1 \\
                      \vdots \\
                      x_n \end{array} \right) =
(f_n \circ \rho_{n+1,n})(x_1, \dots ,x_{n+1})
\]
which concludes the proof. \hfill $\Box$

\smallskip
We are now in a position to show that the exact analogue of 
the classical Floquet-Liapunov theorem is true in
$\mbox{{\Bbb C}}^{\infty}$. More precisely, we have 

\begin{thm}  
Consider the equation (1) with a periodic coefficient of the form
$A : I \longrightarrow {\cal L}(\mbox{{\Bbb C}}^{\infty})$. Then,
the (equivalent) conditions of Theorem 2.3 are always satisfied. 
\end{thm}

\noindent
{\bf Proof.} Clearly, it suffices to prove the existence of a
projective system
$(B_i)_{i \in \mbox{{\Bbs N}}} \in H(\mbox{{\Bbb C}}^{\infty})$
such that $\mbox{Exp}((B_i)_{i \in \mbox{{\Bbs N}}})=\alpha^{\#}(1)$.
As we have seen in Section 2 (see also (3)), $\alpha^{\#}(1) =
\displaystyle{\lim_{\longleftarrow}}(\alpha_n^{\#}(1))$. For the sake
of simplicity we set
\[
M_n := \alpha_n^{\#}(1)  \in  GL(n,\mbox{{\Bbb C}}).
\]
The previous matrices satisfy condition (ii) of Lemma 3.1. We 
define 
\[
B_1 := \log M_1 = \log\lambda_1  \in  \mbox{{\Bbb C}},
\]
where $\log\lambda_1$ is arbitrarily chosen but fixed. By induction,
we set
\[
B_{n+1} := \left( \begin{array}{cc}
                  B_n & 0 \\
                  y_n & \log\lambda_n
                  \end{array} \right), \quad  n \in \mbox{{\Bbb N}}.
\]
Here $y_n$  is determined through the following equivalent conditions:
\[
\exp_{n+1}(B_{n+1}) = M_{n+1} \quad \Leftrightarrow
\]
\[
    \Leftrightarrow \quad \left( \begin{array}{cc}
              \exp_n(B_n)_{_{}} & 0_{_{}} \\
              y_n \cdot \displaystyle{\big(\sum_{k \geq 1} \frac{1}{k!}
              (\sum_{j=1}^k} B_n^{k-j}\cdot (\log\lambda_n)^{j-1})\big)
               & \lambda_n
              \end{array} \right) = M_{n+1} \quad \Leftrightarrow
\]
\[
     \Leftrightarrow \quad  y_n \cdot \big(\sum_{k \geq 1} \frac{1}{k!}
              (\sum_{j=1}^k B_n^{k-j}\cdot (\log\lambda_n)^{j-1})\big)
                            = \mu_n.
\]
Note that the last equation can be solved since the ($n \times n$) matrix 
figuring in it is always non singular. Indeed, the above mentioned 
matrix is triangular, since $B_n$ is. Therefore, its determinant 
equals the product of
\[
\gamma_i :=
\sum_{k \geq 1} \frac{1}{k!}
        (\sum_{j=1}^{k} b_i^{k-j} \cdot (\log\lambda_n)^{j-1}));
              \qquad i = 1,\dots n,
\]
where $b_i$ are the diagonals of $B_n$ satisfying $\exp(b_i)=\lambda_i$.
If $\lambda_n \neq \lambda_i$ ($1 \leq i \leq n-1$), then
$\gamma_i = (\lambda_i - \lambda_n)/(b_i-\log\lambda_n) \neq 0$, for
every $i$. If $\lambda_n = \lambda_i$, for some
$i$, then $\gamma_i = \lambda_n$ which is also non zero since
$M_{n+1} \in GL(n+1,\mbox{{\Bbb C}})$. As a
result, the determinant of 
\[
\sum_{k \geq 1} \frac{1}{k!}\cdot 
            (\sum_{j=1}^k B_n^{k-j}\cdot (\log\lambda_n)^{j-1})
\]
is, in any case, non zero. 

Using once again Lemma 3.1, we check that
$(B_n)_{n \in \mbox{{\Bbs N}}} \in H(\mbox{{\Bbb C}}^{\infty}$).
Moreover, 
\[
\mbox{Exp}((B_n)_{n \in \mbox{{\Bbs N}}}) =
(\exp_n(B_n))_{n \in \mbox{{\Bbs N}}} = \alpha^{\#}(1)
\]
and the proof is now complete. \hfill $\Box$

\medskip
\noindent
{\bf Remark.} Clearly, the proof of the above Theorem is essentially 
based on the appropriate choice of the logarithms $B_n$,
$n \in \mbox{{\Bbb N}}$,
satisfying condition (ii) of Theorem 2.3, instead of considering 
arbitrarily chosen logarithms of $\alpha_n^{\#}(1)$.

\section{A geometric generalization} 

In this section we briefly describe an application of the 
previous methods to the geometric framework of Fr\'echet fiber 
bundles. Since the complete details are beyond the scope of this 
paper, we restrict ourselves to a mere outline of these ideas, the 
Banach analogues of which can be found in \cite{13,14}. 

Assume that the coefficient $A$ of (1) is periodic. With the 
notations of Section 2, if we interpret $\mbox{{\Bbb Z}}$ as the
fundamental group of the unit circle; that is, $\mbox{{\Bbb Z}} \cong
\pi_1(S^1))$, and $\mbox{{\Bbb R}}$ as the universal covering
space of $S^1$ (: $\mbox{{\Bbb R}} \cong \widetilde{S^1}$), then the
corresponding fundamental solution of (1), given by 
\[
\Phi  : \mbox{{\Bbb R}} \longrightarrow H^o(\mbox{{\Bbb F}})
: t \mapsto (\Phi_i(t))_{i \in \mbox{{\Bbs N}}},
\]
determines the $H(\mbox{{\Bbb F}})$-valued smooth 1-form 
\[
\widetilde{\theta} := d\Phi \cdot \Phi^{-1}   \in
\Lambda^1(\widetilde{S^1},H(\mbox{{\Bbb F}})).
\]
We recall that $d\Phi \cdot \Phi^{-1}$ denotes the right total differential 
\[
(d\Phi \cdot \Phi^{-1})_x = (dR_{\Phi(x)^{-1}})_{\Phi(x)} \circ
(d\Phi)_x;  \quad x \in \mbox{{\Bbb R}} \cong \widetilde{S^1},
\]
where $R_{\Phi(x)^{-1}}$ is the right translation of
$H^o(\mbox{{\Bbb F}})$ by $\Phi(x)^{-1}$.

Using projective limits, along with \cite{13} and \cite{14}, we check 
that there exists an {\em integrable\/} $\theta \in
\Lambda^1(S^1,H(\mbox{{\Bbb F}}))$ (that is,
$d\theta = 1/2 \cdot [\theta,\theta])$ such that
\[
p^*\theta = \widetilde{\theta},
\]
where $p : \widetilde{S^1} \cong \mbox{{\Bbb R}} \longrightarrow S^1 :
t \mapsto \exp(2\pi it)$. Then, we can find a unique
connection $\omega_{\Phi}$ on the trivial bundle
$\ell_o = (S^1 \times H^o(\mbox{{\Bbb F}}),
H^o(\mbox{{\Bbb F}}),S^1,pr_1)$
with unique (local) connection form $-\,\theta$, i.e. $-\,\theta =
\sigma^*\omega_{\Phi}$, if $\sigma : S^1 \longrightarrow
S^1 \times H^o(\mbox{{\Bbb F}})$ is the natural section
of the bundle. It turns out that $\omega_{\Phi}$ is a flat
connection with corresponding holonomy homomorphism coinciding, up 
to conjugation, with $\alpha^{\#}$. Moreover, the coefficient of
(1) is related with $\omega$ by 
\[
A(t) = ((\sigma \circ p)^* \omega_{\Phi})(\partial_t),
\]
if $\partial$ is the fundamental vector field of $\mbox{{\Bbb R}}$.
The converse is also true. Thus, there is a bijection between flat
connections on $\ell_o$ and equations of type (1) with periodic
coefficients. Similarly, we can construct an $\omega_F$ for any
morphism $F : \mbox{{\Bbb R}} \longrightarrow
H^o(\mbox{{\Bbb F}})$. Now, the
geometric analogue of Theorem 2.3 can be stated as follows: 

\smallskip
\begin{flushright}
\begin{minipage}{11cm}
{\em If the monodromy homomorphism $\alpha^{\#}$ of (1) (with periodic 
coefficient) extends to a morphism $F : \mbox{{\Bbb R}}
\longrightarrow H^o(\mbox{{\Bbb F}})$, then  the
connection $\omega_F$ corresponds to a constant coefficient $B$.} 
\end{minipage}
\end{flushright}

\smallskip
\noindent
We note that $\omega_F$, $\omega_{\Phi}$ are gauge-equivalent and
that the coefficient $B$ is precisely the coefficient of Theorem 2.3. 

The previous result could motivate an analogous study on 
arbitrary (non trivial) $H^o(\mbox{{\Bbb F}})$-bundles. This idea
has been applied in \cite{15}

\medskip
\flushright{%
\begin{minipage}{5cm}
University of Athens \\
Department of Mathematics \\
Panepistimiopolis \\
Athens 157 84, Greece\\
ggalanis@atlas.uoa.gr \\
evassil@atlas.uoa.gr
\end{minipage}}


\begin{thebibliography}{99}

\bibitem{1}
{\sc H.~Amann}, ``Ordinary differential equations'', W. de Gruyter,
   Berlin, 1990. 

\bibitem{2}
{\sc V.I.~Bogachev}, {\em Deterministic and stochastic differential
equations in infinite-dimensional spaces\/}, Acta Applicandae Mathematicae
{\bf 40} (1995), 25--93.

\bibitem{3}
{\sc H.~Cartan}, ``Differential Calculus'', Hermann, Paris, 1971. 

\bibitem{4}
{\sc G.~Galanis}, {\em On a type of linear differential equations in 
Fr\'echet spaces\/}, Ann. Scuola Norm. Sup. Pisa, Ser. 4,
{\bf 24} (1997), 501--510

\bibitem{5}
{\sc R.S.~Hamilton}, {\em The Inverse function theorem of Nash and Moser\/}, 
Bull. Amer. Math. Soc. {\bf 7} (1982), 65--222. 

\bibitem{6}
{\sc S.G.~Krein -- N.I.~Jackin}, ``Linear differential equations on
manifolds'', Veronez. Gos. Univ., 1980 (in Russian). 

\bibitem{7}
{\sc R.~Lemmert}, {\em On ordinary differential equations in locally convex 
spaces\/}, Nonlinear Analysis, Theory, Methods and Applications {\bf 10}
(1986), 1385--1390. 

\bibitem{8}
{\sc J.A.~Leslie}, {\em On a differential structure for the group of
diffeomorphisms\/}, Topology {\bf 6} (1967), 263--271.

\bibitem{9}
{\sc J.L.~Massera -- J.J.Schaeffer}, ``Linear differential equations 
and function spaces'', Academic Press, London, 1966. 

\bibitem{10}
{\sc L.S.~Pontryagin}, ``Ordinary differential equations'', Addison-Wesley,
Reading Mass., 1962.

\bibitem{11}
{\sc H.H.~Schaeffer}, ``Topological Vector Spaces'', Springer-Verlag, 
1980. 

\bibitem{12}
{\sc J.J.~Schaeffer}, {\em On Floquet's theorem in Hilbert spaces\/}, Bull. 
Amer. Math. Soc. {\bf 70} (1964), 243-245. 

\bibitem{13}
{\sc E.~Vassiliou}, {\em Floquet-type connections on principal bundles\/}, 
Tensor N.S. {\bf 43} (1986), 189--195. 

\bibitem{14}
{\sc E.~Vassiliou}, {\em Total differential equations and the structure of 
fibre bundles\/}, Bull. Greek. Math. Soc. {\bf 27} (1986), 149-159. 

\bibitem{15}
{\sc E.~Vassiliou -- G.Galanis}, {\em Flat Fr\'echet principal bundles and 
their holonomy homomorphisms\/} (submitted). 

\bibitem{16}
{\sc H.K.~Wilson}, ``Ordinary differential equations'', Addison-Wesley, 
Reading Mass., 1971. 


\end{thebibliography}
\end{document}